\documentclass[11pt]{amsart}
\usepackage{amssymb}
\usepackage{epsfig}
\newtheorem{theorem}{Theorem}

\newtheorem{proposition}[theorem]{Proposition}
\newtheorem{corollary}[theorem]{Corollary}
\newtheorem{definition}[theorem]{Definition}
\newtheorem{conj}[theorem]{Conjecture}
\newtheorem{question}[theorem]{Question}

\newcommand{\C}{\mathbb{C}}
\newcommand{\Z}{\mathbb{Z}}
\newcommand{\R}{\mathbb{R}}
\newcommand{\CP}{\mathbb{CP}}
\title[Symplectic 4-manifolds, plane curves, and isotopy]%
{Symplectic 4-manifolds, singular plane curves, and isotopy problems}
\author{Denis Auroux}
\address{Department of Mathematics, M.I.T., Cambridge MA 02139, USA}
\email{auroux@math.mit.edu}
\thanks{From lectures given at the 2004 Clay Mathematics Institute Summer
School on Floer Homology, Gauge Theory, and Low Dimensional Topology at
the Alfr\'ed R\'enyi Institute, Budapest, Hungary
({\tt www.claymath.org/programs/summer\underline{\ }school/2004/}).\\
\indent Partially supported by NSF grant DMS-0244844.}

\begin{document}
\begin{abstract}
We give an overview of various recent results concerning the topology of
symplectic 4-manifolds and singular plane curves, using branched covers
and isotopy problems as a unifying theme. While this paper does not
contain any new results, we hope that it can serve as an introduction
to the subject, and will stimulate interest in some of the open
questions mentioned in the final section.
\end{abstract}

\maketitle

\section{Introduction}

An important problem in 4-manifold topology is to understand which manifolds
carry symplectic structures (i.e., closed non-degenerate 2-forms), and to
develop invariants that can distinguish symplectic manifolds. Additionally,
one would like to understand to what extent the category of symplectic
manifolds is richer than that of K\"ahler (or complex projective) manifolds.
Similar questions may be asked about singular curves inside, e.g., the
complex projective plane. The two types of questions are related to each
other via symplectic branched covers.
 
A branched cover of a symplectic 4-manifold with a (possibly singular)
symplectic branch curve carries a natural symplectic structure. Conversely,
using approximately holomorphic techniques it can be shown that every
compact symplectic 4-manifold is a branched cover of the complex projective
plane, with a branch curve presenting nodes (of both orientations) and
complex cusps as its only singularities (cf.\ \S \ref{sec:covers}).
The topology of the 4-manifold
and that of the branch curve are closely related to each other; for
example, using braid monodromy techniques to study the branch curve, one
can reduce the classification of symplectic 4-manifolds to a (hard) question
about factorizations in the braid group (cf.\ \S \ref{sec:bmf}).
Conversely, in some examples the topology of the branch curve complement
(in particular its fundamental group) admits a simple description in terms
of the total space of the covering (cf.\ \S \ref{sec:pi1}).
                                                                                
In the language of branch curves, the failure of most symplectic manifolds
to admit integrable complex structures translates into the failure of most
symplectic branch curves to be isotopic to complex curves. While the
symplectic isotopy problem has a negative answer for plane curves with
cusp and node singularities, it is interesting to investigate this failure
more precisely. Various partial results have been obtained recently about
situations where isotopy holds (for smooth curves; for curves of low
degree), and about isotopy up to stabilization or regular homotopy
(cf.\ \S \ref{sec:isotopy}). On the
other hand, many known examples of non-isotopic curves can be understood
in terms of twisting along Lagrangian annuli (or equivalently, Luttinger
surgery of the branched covers), leading to some intriguing open questions
about the topology of symplectic 4-manifolds versus that of K\"ahler surfaces.

\section{Background}

In this section we review various classical facts about symplectic
manifolds; the reader unfamiliar with the subject is referred to the
book \cite{McS} for a systematic treatment of the material.

Recall that a {\it symplectic form} on a smooth manifold is a 2-form
$\omega$ such that $d\omega=0$ and $\omega\wedge\dots\wedge \omega$ is a
volume form. The prototype of a symplectic form is the 2-form $\omega_0=\sum
dx_i\wedge dy_i$ on $\R^{2n}$. In fact, one of the most classical results
in symplectic topology, Darboux's theorem, asserts that every symplectic
manifold is locally symplectomorphic to $(\R^{2n},\omega_0)$: hence, unlike
Riemannian metrics, symplectic structures have no local invariants.

Since we are interested primarily in compact examples, let us mention 
compact oriented surfaces (taking $\omega$ to be an arbitrary area form),
and the complex projective space $\CP^n$ (equipped with the Fubini-Study
K\"ahler form). More generally, since any submanifold to which $\omega$
restricts non-degenerately inherits a symplectic structure, all complex
projective manifolds are symplectic. However, the symplectic category
is strictly larger than the complex projective category, as first
evidenced by Thurston in 1976 \cite{Th}. In 1994 Gompf obtained the
following spectacular result using the {\it symplectic sum} construction
\cite{Go1}:

\begin{theorem}[Gompf]
Given any finitely presented group $G$, there exists a compact symplectic
4-manifold $(X,\omega)$ such that $\pi_1(X)\simeq G$.
\end{theorem}

Hence, a general symplectic manifold cannot be expected to carry a complex
structure; however, we can equip it with a compatible {\it almost-complex}
structure, i.e.\
there exists $J\in\mathrm{End}(TX)$ such that $J^2=-\mathrm{Id}$ and
$g(\cdot,\cdot):=\omega(\cdot,J\cdot)$ is a Riemannian metric. Hence, at
any given point $x\in X$ the tangent space $(T_xX,\omega,J)$ can be
identified with $(\C^n,\omega_0,i)$, but there is no control over the
manner in which $J$ varies from one point to another ($J$ is not {\it
integrable}). In particular, the $\bar\partial$ operator associated to $J$
does not satisfy $\bar\partial^2=0$, and hence there are no local
holomorphic coordinates.
\medskip

An important problem in 4-manifold topology is to understand the hierarchy
formed by the three main classes of compact oriented 4-manifolds: (1)
complex projective, (2) symplectic, and (3) smooth. Each class is a proper
subset of the next one, and many obstructions and examples are known,
but we are still very far from understanding
what exactly causes a smooth 4-manifold to admit a symplectic structure, or
a symplectic 4-manifold to admit an integrable complex structure.

One of the main motivations to study symplectic 4-manifolds
is that they retain some (but not all) features of complex projective
manifolds: for example the structure of their Seiberg-Witten invariants,
which in both cases are non-zero and count certain embedded curves
\cite{Ta1,Ta2}. At the same time, every compact oriented smooth 4-manifold with
$b_2^+\ge 1$ admits a ``near-symplectic'' structure, i.e.\ a closed 2-form
which vanishes along a union of circles and is symplectic over the
complement of its zero set \cite{GK,Ho1}; and it appears that some
structural properties of symplectic manifolds carry over to the world of
smooth 4-manifolds (see e.g.\ \cite{Ta3,Asinglp}).

Many new developments have contributed to improve our understanding of
symplectic 4-manifolds over the past ten years (while results are much
scarcer in higher dimensions). Perhaps the most important source of
new results has been the study of pseudo-holomorphic curves in their various
incarnations: Gromov-Witten invariants, Floer homology, \dots
(for an overview of the subject see \cite{McS2}). At the same time,
gauge theory (mostly Seiberg-Witten theory, but also more recently
Ozsvath-Szabo theory) has made it possible to identify various {\it
obstructions} to the existence of symplectic structures in dimension 4
(cf.\ e.g.\ \cite{Ta1,Ta2}). On the other hand, various new constructions,
such as link surgery \cite{FS1}, symplectic sum \cite{Go1}, and symplectic
rational blowdown \cite{Sy} have made it possible to exhibit interesting
families of non-K\"ahler symplectic 4-manifolds. In a slightly different
direction, approximately holomorphic geometry (first introduced by Donaldson
in \cite{Do1}) has made it
possible to obtain various structure results, showing that symplectic
4-manifolds can be realized as symplectic Lefschetz pencils \cite{Do2} or
as branched covers of $\CP^2$ \cite{Au2}. In the rest of this paper we
will focus on this latter approach, and discuss the topology of {\it
symplectic branched covers} in dimension 4.

\section{Symplectic branched covers}\label{sec:covers}

Let $X$ and $Y$ be compact oriented 4-manifolds, and assume that $Y$ carries
a symplectic form $\omega_Y$.

\begin{definition}
A smooth map $f:X\to Y$ is a {\em symplectic branched covering} if given any
point $p\in X$ there exist neighborhoods $U\ni p$, $V\ni f(p)$, and local
coordinate charts $\phi:U\to\C^2$
$($orientation-preserving$)$ and $\psi:V\to\C^2$
$($adapted to $\omega_Y$, i.e.\ such that $\omega_Y$ restricts positively
to any complex line in $\C^2)$, in which $f$ is given by one of:
\smallskip

$(i)$ $(x,y)\mapsto (x,y)$ $($local diffeomorphism$)$,

$(ii)$ $(x,y)\mapsto (x^2,y)$ $($simple branching$)$,

$(iii)$ $(x,y)\mapsto (x^3-xy,y)$ $($ordinary cusp$)$.
\end{definition}

These local models are the same as for the singularities of a generic 
holomorphic map from $\C^2$ to itself, except that the requirements on the
local coordinate charts have been substantially weakened.
The {\it ramification curve} $R=\{p\in X,\ \det(df)=0\}$ is a smooth
submanifold of $X$, and its image $D=f(R)$ is the {\it branch curve}, described
in the local models by the equations $z_1=0$ for $(x,y)\mapsto (x^2,y)$
and $27z_1^2=4z_2^3$ for $(x,y)\mapsto (x^3-xy,y)$. The conditions imposed
on the local coordinate charts imply that $D$ is a
{\it symplectic curve} in $Y$ (i.e., $\omega_{Y|TD}>0$ at every point of
$D$). Moreover the restriction of $f$ to $R$
is an immersion everywhere except at the cusps. Hence, besides the ordinary
complex cusps imposed by the local model, the only generic singularities
of $D$ are transverse double points (``nodes''), which may occur with either the complex
orientation or the anti-complex orientation.

We have the following result \cite{Au2}:

\begin{proposition}\label{prop:au2}
Given a symplectic branched covering $f:X\to Y$, the manifold $X$ inherits
a natural symplectic structure $\omega_X$, canonical up to isotopy, in the
cohomology class $[\omega_X]=f^*[\omega_Y]$.
\end{proposition}

The symplectic form $\omega_X$ is constructed by adding to $f^*\omega_Y$
a small multiple of an exact form $\alpha$ with the property that, at
every point of $R$, the restriction of $\alpha$ to $\mathrm{Ker}(df)$ is
positive. Uniqueness up to isotopy follows from the
convexity of the space of such exact 2-forms and Moser's theorem.

Conversely, we can realize every compact symplectic 4-manifold
as a symplectic branched cover of $\CP^2$ \cite{Au2}, at least if we
assume {\it integrality}, i.e.\ if we require that
$[\omega]\in H^2(X,\Z)$, which does not place any additional restrictions
on the diffeomorphism type of $X$:

\begin{theorem}\label{thm:au2}
Given an integral compact symplectic 4-manifold $(X^4,\omega)$ and an
integer $k\gg 0$, there exists a symplectic branched covering
$f_k:X\to\CP^2$, canonical up to isotopy if $k$ is sufficiently large.
\end{theorem}

Moreover, the natural symplectic structure induced on $X$ by the
Fubini-Study K\"ahler form and $f_k$ (as given
by Proposition \ref{prop:au2}) agrees with $\omega$ up to isotopy and
scaling (multiplication by~$k$).

The main tool in the construction of the maps $f_k$ is {\it approximately
holomorphic geometry} \cite{Do1,Do2,Au2}. Equip
$X$ with a compatible almost-complex structure, and consider a complex
line bundle $L\to X$ such that $c_1(L)=[\omega]$: then for $k\gg 0$ the 
line bundle $L^{\otimes k}$ admits many approximately holomorphic sections,
i.e.\ sections such that $\sup |\bar\partial s|\ll\sup |\partial s|$.
Generically, a triple of such sections $(s_0,s_1,s_2)$ has no common zeroes,
and determines a projective map $f:p\mapsto [s_0(p)\!:\!s_1(p)\!:\!s_2(p)]$.
Theorem \ref{thm:au2} is then proved by constructing triples of sections
which satisfy suitable transversality estimates, ensuring that the structure
of $f$ near its critical locus is the expected one \cite{Au2}. (In the complex
case it would be enough to pick three generic holomorphic sections,
but in the approximately holomorphic context one needs to work harder and
obtain uniform transversality estimates on the derivatives of $f$.)

Because for large $k$ the maps $f_k$ are canonical up to isotopy through
symplectic branched covers, the topology of $f_k$ and of its branch curve
$D_k$ can be used to define invariants of the symplectic
manifold $(X,\omega)$. The only generic singularities of the plane curve $D_k$
are nodes (transverse double points) of either orientation and complex
cusps, but in a generic one-parameter family of branched covers pairs of nodes
with opposite orientations may be cancelled or created. However, recalling
that a node of $D_k$ corresponds to the occurrence of two simple branch
points in a same fiber of $f_k$, the creation of a pair of nodes can only
occcur in a manner compatible with the branched covering structure, i.e.\
involving disjoint sheets of the covering. Hence, for large $k$ the sequence
of branch curves $D_k$ is, up to isotopy (equisingular deformation among
symplectic curves), cancellations and admissible creations of pairs of
nodes, an invariant of $(X,\omega)$.

The ramification curve of $f_k$ is just a smooth connected symplectic curve
representing the homology class Poincar\'e dual to $3k[\omega]-c_1(TX)$,
but the branch curve $D_k$ becomes more and more complicated as $k$
increases: in terms of the symplectic volume and Chern numbers of $X$,
its degree (or homology class) $d_k$,
genus $g_k$, and number of cusps $\kappa_k$ are given by
$$d_k=3k^2\,[\omega]^2-k\,c_1\cdot [\omega],\qquad
2g_k-2=9 k^2\,[\omega]^2-9 k\,c_1\cdot [\omega]+2c_1^2,$$
$$\kappa_k=12k^2\,[\omega]^2-9k\,c_1\cdot [\omega]+2c_1^2-c_2.$$
It is also worth mentioning that, to this date, there is no evidence suggesting
that negative nodes actually do occur in these high degree branch curves;
our inability to rule our their presence might well be a shortcoming of the
approximately holomorphic techniques, rather than an intrinsic feature of
symplectic 4-manifolds. So in the following sections we will occasionally
consider the more conventional problem of understanding isotopy classes
of curves presenting only positive nodes and cusps, although most of the
discussion applies equally well to curves with negative nodes.
\medskip

Assuming that the topology of the branch curve is understood (we will
discuss how to achieve this in the next section), one still needs to
consider the branched covering $f$ itself. The structure of $f$ is
determined by its {\it monodromy morphism} $\theta:\pi_1(\CP^2-D)\to S_N$,
where $N$ is the degree of the covering $f$. Fixing a base point $p_0\in
\CP^2-D$, the image by $\theta$ of a loop $\gamma$ in the complement of $D$
is the permutation of the fiber $f^{-1}(p_0)$ induced by the monodromy of $f$
along $\gamma$. (Since viewing this permutation as an element of $S_N$
depends on the choice of an identification
between $f^{-1}(p_0)$ and $\{1,\dots,N\}$, the morphism $\theta$ is only
well-defined up to conjugation by an element of $S_N$.) By Proposition
\ref{prop:au2}, the isotopy class of the branch curve $D$ and the monodromy
morphism $\theta$ determine completely the symplectic 4-manifold $(X,\omega)$
up to symplectomorphism.

Consider a loop $\gamma$ which bounds a small topological disc
intersecting $D$ transversely once: such a loop plays a role similar to
the meridian of a knot, and is called a {\it geometric
generator} of $\pi_1(\CP^2-D)$. Then $\theta(\gamma)$ is a transposition (because
of the local model near a simple branch point). Since the image of $\theta$
is generated by transpositions and acts transitively on the fiber (assuming
$X$ to be connected), $\theta$ is a surjective group homomorphism. Moreover,
the smoothness of $X$ above the singular points of $D$ imposes certain
compatibility conditions on $\theta$. Therefore, not every singular plane
curve can be the branch curve of a smooth covering; moreover, the morphism
$\theta$, if it exists, is often unique (up to conjugation in $S_N$).
In the case of algebraic curves, this uniqueness property, which holds
except for a finite list of well-known counterexamples, is known as
Chisini's conjecture, and was essentially proved by Kulikov a few years
ago \cite{Ku}.

The upshot of the above discussion is that, in order to understand symplectic
4-manifolds, it is in principle enough to understand singular plane curves.
Moreover, if the branch curve of a symplectic covering $f:X\to \CP^2$ happens
to be a complex curve, then the integrable complex structure of $\CP^2$ can be
lifted to an integrable complex structure on $X$, compatible with the
symplectic structure; this implies
that $X$ is a complex projective surface. So, considering the branched
coverings constructed in Theorem \ref{thm:au2}, we have:

\begin{corollary}\label{cor:au2}
For $k\gg 0$ the branch curve $D_k\subset\CP^2$ is isotopic to a complex
curve (up to node cancellations) if and only if $X$ is a complex projective
surface.
\end{corollary}

This motivates the study of the {\it symplectic isotopy problem}, which
we will discuss in \S \ref{sec:isotopy}. For now we focus on the use
of braid monodromy invariants to study the topology of singular plane
curves. In the present context, the goal of this approach is to reduce the
classification
of symplectic 4-manifolds to a purely algebraic problem, in a manner
vaguely reminiscent of the role played by Kirby calculus in the
classification of smooth 4-manifolds; as we shall see below, 
representing symplectic 4-manifolds as branched covers of $\CP^2$
naturally leads one to study the
calculus of factorizations in braid groups.

\section{The topology of singular plane curves}\label{sec:bmf}

The topology of singular algebraic plane curves has been studied
extensively since Zariski. One of the main tools is the notion of
{\it braid monodromy} of a plane curve, which has been used in particular
by Moishezon and Teicher in many papers since the early 1980s in order
to study branch curves of generic projections of complex projective surfaces
(see \cite{Te1} for a detailed overview). 
Braid monodromy techniques can be applied to the more general case of
{\it Hurwitz curves} in ruled surfaces, i.e.\ curves which behave
in a generic manner with respect to the ruling. In the case of $\CP^2$,
we consider the projection $\pi:\CP^2-\{(0:0:1)\}\to \CP^1$ given by
$(x:y:z)\mapsto (x:y)$.

\begin{definition}\label{def:hurwitz}
 A curve $D\subset\CP^2$ $($not passing through $(0\!:\!0\!:\!1))$ is a
Hurwitz curve (or braided curve) if $D$ is positively transverse to the fibers of $\pi$
everywhere except at finitely many points where $D$ is smooth and
non-degenerately tangent to the fibers.
\end{definition}

\begin{figure}[t]
\begin{center}
\setlength{\unitlength}{0.7mm}
\begin{picture}(80,52)(-40,-12)
\put(0,-2){\vector(0,-1){8}}
\put(2,-6){$\pi:(x:y:z)\mapsto (x:y)$}
\put(-40,-15){\line(1,0){80}}
\put(-38,-12){$\CP^1$}
\put(-40,0){\line(1,0){80}}
\put(-40,40){\line(1,0){80}}
\put(-40,0){\line(0,1){40}}
\put(40,0){\line(0,1){40}}
\put(-38,33){$\CP^2-\{0\!:\!0\!:\!1\}$}
\put(27,31){$D$}
\multiput(-20,20)(0,-2){18}{\line(0,-1){1}}
\multiput(-5,20)(0,-2){18}{\line(0,-1){1}}
\multiput(15,15)(0,-2){9}{\line(0,-1){1}}
\multiput(15,-9)(0,-2){3}{\line(0,-1){1}}
\put(-20,-15){\circle*{1}}
\put(-5,-15){\circle*{1}}
\put(15,-15){\circle*{1}}
\qbezier[140](25,35)(5,30)(-5,20)
\qbezier[60](-5,20)(-10,15)(-15,15)
\qbezier[60](-15,15)(-20,15)(-20,20)
\qbezier[60](-20,20)(-20,25)(-15,25)
\qbezier[60](-15,25)(-10,25)(-5,20)
\qbezier[100](-5,20)(0,15)(15,15)
\qbezier[250](15,15)(5,15)(-30,5)
\put(-20,20){\circle*{1}}
\put(-5,20){\circle*{1}}
\put(15,15){\circle*{1}}
\end{picture}
\end{center}
\caption{A Hurwitz curve in $\CP^2$}
\end{figure}
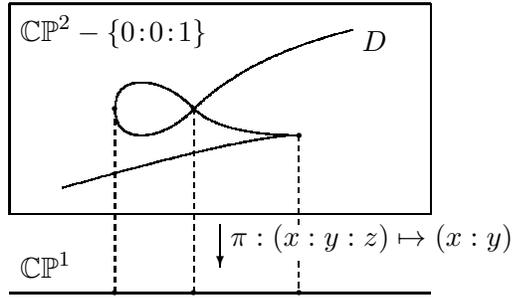

The projection $\pi$ makes $D$ a singular branched cover of $\CP^1$,
of degree $d=\deg D=[D]\cdot[\CP^1]$. Each fiber of $\pi$ is a
complex line $\ell\simeq \C\subset\CP^2$, and if $\ell$ does not pass
through any of the singular points of $D$ nor any of its vertical
tangencies, then $\ell\cap D$ consists of $d$ distinct points.
We can trivialize the fibration $\pi$ over an affine subset
$\C\subset\CP^1$, and define
the {\it braid monodromy morphism}
$$\rho:\pi_1(\C-\mathrm{crit}(\pi_{|D}))\to B_d.$$
Here $B_d$ is the Artin braid group on $d$ strings (the fundamental group
of the configuration space $\mathrm{Conf}_d(\C)$
of $d$ distinct points in $\C$), and for any loop
$\gamma$ the braid $\rho(\gamma)$ describes the
motion of the $d$ points of $\ell\cap D$ inside the
fibers of $\pi$ as one moves along the loop $\gamma$.

Equivalently, choosing an ordered system of arcs generating the free group
$\pi_1(\C-\mathrm{crit}(\pi_{|D}))$, one can express the braid monodromy
of $D$ by a {\it factorization} $$\Delta^2=\prod_{i} \rho_i$$ of the central
element $\Delta^2$ (representing a full rotation by $2\pi$) in $B_d$, where each factor $\rho_i$ is the monodromy
around one of the special points (cusps, nodes, tangencies) of $D$.

A same Hurwitz curve can be described by different factorizations of
$\Delta^2$ in $B_d$: switching to a different ordered system of generators
of $\pi_1(\C-\mathrm{crit}(\pi_{|D}))$ affects the collection of factors
$\langle \rho_1,\dots,\rho_r\rangle $ by a sequence of {\it Hurwitz moves},
i.e.\ operations of the form 
$$\langle \rho_1,\,\cdots,\rho_i,\rho_{i+1},\,\cdots,\rho_r\rangle \,
\longleftrightarrow\, \langle \rho_1,\,\cdots,(\rho_i\rho_{i+1}\rho_i^{-1}),
\rho_i,\,\cdots,\rho_r\rangle;
$$
and changing the identification between the reference fiber
$(\ell,\ell\cap D)$ of $\pi$ and the base point in $\mathrm{Conf}_d(\C)$
affects braid monodromy by a {\it global conjugation}
$$\langle\rho_1,\,\cdots,\rho_r\rangle \,\longleftrightarrow\,
\langle b^{-1}\rho_1 b,\,\cdots,b^{-1}\rho_r b\rangle.
$$
For Hurwitz curves whose only singularities are cusps and nodes (of either
orientation), or more generally curves with $A_n$ (and $\overline{A}_n$)
singularities, the braid monodromy factorization determines the isotopy
type completely (see for example \cite{KK}). Hence,
determining whether two given Hurwitz curves are isotopic among Hurwitz
curves is equivalent to determining whether two given factorizations of
$\Delta^2$ coincide up to Hurwitz moves and global conjugation. 
\medskip

It is easy to see that any Hurwitz curve in $\CP^2$ can be made symplectic
by an isotopy through Hurwitz curves: namely, the image of any Hurwitz curve
by the rescaling map $(x:y:z)\mapsto (x:y:\lambda z)$ is a Hurwitz curve,
and symplectic for $|\lambda|\ll 1$. On the other hand, a refinement
of Theorem \ref{thm:au2} makes it possible to assume without loss of
generality that the branch curves
$D_k\subset\CP^2$ are Hurwitz curves \cite{AK}. So, from now on we
can specifically consider symplectic coverings with Hurwitz branch
curves. In this setting,
braid monodromy gives a purely combinatorial description of 
the topology of compact (integral) symplectic 4-manifolds.

The braid monodromy of the branch curves $D_k$ given by Theorem
\ref{thm:au2} can be computed explicitly for various families of complex
projective surfaces (non-K\"ahler examples are currently beyond reach).
In fact, in  the complex case the branched coverings $f_k$ are
isotopic to generic projections of projective embeddings. Accordingly, most
of these computations rely purely on methods from algebraic geometry, using
the degeneration techniques extensively developed by Moishezon and Teicher
(see \cite{AGTV,Mo1,Mo2,MRT,Ro,Te1,Te2} and references within); but
approximately holomorphic methods can be used to simplify the calculations
and bring a whole new range of examples
within reach \cite{ADKY}. This includes some complex surfaces of general
type which are mutually homeomorphic and have identical Seiberg-Witten
invariants but of which it is unknown whether they are symplectomorphic
or even diffeomorphic (the {\it Horikawa surfaces}).

However, the main obstacle standing in the way of this approach to the
topology of symplectic 4-manifolds is the intractability of the so-called
``Hurwitz problem'' for braid monodromy factorizations: namely, there is
no algorithm to decide whether two given braid monodromy factorizations
are identical up to Hurwitz moves. Therefore, since we are unable to compare
braid monodromy factorizations, we have to extract the information contained
in them by indirect means, via the introduction of more manageable (but less
powerful) invariants.

\section{Fundamental groups of branch curve complements}\label{sec:pi1}

The idea of studying algebraic plane curves by determining the fundamental
groups of their complements is a very classical one, which goes back
to Zariski and Van Kampen. More recently, Moishezon and Teicher have
shown that fundamental groups of branch curve complements can be used as
a major tool to further our understanding of complex projective surfaces
(cf.\ e.g.\ \cite{Mo1,MT,Te1}). By analogy with the situation for knots
in $S^3$, one expects the topology of the complement to carry a lot of
information about the curve; however in this case the fundamental group
does not determine the isotopy type. For an algebraic curve in $\CP^2$, or
more generally for a Hurwitz curve, the fundamental group of the complement
is determined in an explicit manner by the braid monodromy factorization,
via the Zariski-Van Kampen theorem. Hence, calculations of fundamental
groups of complements usually rely on braid monodromy techniques.

A close examination of the available data suggests that, contrarily to what
has often been claimed, in the specific case of generic projections of
complex surfaces projectively embedded by sections of a sufficiently ample
linear system (i.e.\ taking $k\gg 0$ in Theorem \ref{thm:au2}), the
fundamental group of the branch curve complement may be determined in an
elementary manner by the topology of the surface (see below).

In the symplectic setting, the fundamental group of the complement of the
branch curve $D$ of a covering $f:X\to\CP^2$ is affected by node creation
or cancellation operations. Indeed, adding pairs of nodes (in a manner 
compatible with the monodromy morphism $\theta:\pi_1(\CP^2-D)\to S_N$)
introduces additional commutation relations between geometric generators
of the fundamental group. Hence, it is necessary to consider a suitable
``symplectic stabilization'' of $\pi_1(\CP^2-D)$ \cite{ADKY}:

\begin{definition}\label{def:stabgp}
Let $K$ be the normal
subgroup of $\pi_1(\CP^2-D)$ generated by the commutators $[\gamma,\gamma']$
for all pairs $\gamma,\gamma'$ of geometric generators such that
$\theta(\gamma)$ and $\theta(\gamma')$ are disjoint commuting
transpositions. Then the symplectic stabilization of $\pi_1(\CP^2-D)$ is 
the quotient $\bar{G}=\pi_1(\CP^2-D)/K$.
\end{definition}

Considering the branch curves $D_k$ of the coverings given by Theorem
\ref{thm:au2}, we have the following result \cite{ADKY}:

\begin{theorem}[A.-Donaldson-Katzarkov-Yotov]
For $k\gg 0$, the stabilized group $\bar{G}_k(X,\omega)=
\pi_1(\CP^2-D_k)/K_k$ is an invariant of the symplectic manifold $(X^4,\omega)$.
\end{theorem}

The fundamental group of the complement of a plane branch curve
$D\subset\CP^2$ comes naturally equipped with two morphisms:
the symmetric group valued monodromy homomorphism $\theta$ discussed
above, and the abelianization map $\delta:\pi_1(\CP^2\!-\!D)\to
H_1(\CP^2\!-\!D,\Z)$. Since we only consider irreducible branch curves, we have
$H_1(\CP^2\!-\!D,\Z)\simeq \Z_d$, where $d=\deg D$, and $\delta$ counts the
linking number (mod $d$) with the curve $D$.  The morphisms $\theta$ and
$\delta$ are surjective, but the image of 
$(\theta,\delta):\pi_1(\CP^2-D)\to S_N\times \Z_d$ is the index 2 subgroup
consisting of all pairs $(\sigma,p)$ such that the permutation $\sigma$
and the integer $p$ have the same parity (note that $d$ is always even).
The subgroup $K$ introduced in Definition \ref{def:stabgp} lies in the
kernel of $(\theta,\delta)$; therefore, setting
$G^0=\mathrm{Ker}(\theta,\delta)/K$, we have
an exact sequence
$$1\longrightarrow G^0\longrightarrow \bar{G}\stackrel{(\theta,\delta)}
{\longrightarrow}S_N\times \Z_d\longrightarrow \Z_2\longrightarrow 1.$$
Moreover, assume that the symplectic 4-manifold $X$ is simply connected,
and denote by $L=f^*[\CP^1]$ the pullback of the hyperplane class
and by $K_X=-c_1(TX)$ the canonical class. Then we have the following
result \cite{ADKY}:
\begin{theorem}[A.-Donaldson-Katzarkov-Yotov]\label{thm:adky}
If $\pi_1(X)=1$ then there is a natural surjective homomorphism
$\phi:\mathrm{Ab}(G^0)\twoheadrightarrow (\Z^2/\Lambda)^{N-1}$, where
$\Lambda=\{(L\cdot C, K_X\cdot C),\ C\in H_2(X,\Z)\}\subset\Z^2$.
\end{theorem}

The fundamental groups of the branch curve complements have been computed
for generic polynomial maps to $\CP^2$ on various algebraic surfaces,
using braid monodromy techniques (cf.\ \S \ref{sec:bmf}) and the
Zariski-Van Kampen theorem. Since in the symplectic setting
Theorem \ref{thm:au2} gives uniqueness
up to isotopy only for $k\gg 0$, we restrict ourselves to those examples
for which the fundamental groups have been computed for $\CP^2$-valued maps
of arbitrarily large degree. 

The first such calculations were carried out by Moishezon and Teicher,
for $\CP^2$, $\CP^1\times\CP^1$ \cite{Mo2}, and Hirzebruch surfaces
(\cite{MRT}, see also \cite{ADKY}); the answer is also known for some
specific linear systems on rational surfaces and K3 surfaces realized
as complete intersections (by work of Robb \cite{Ro}, see also related
papers by Teicher et al).
Additionally, the symplectic stabilizations of the fundamental groups
have been computed for all double covers
of $\CP^1\times\CP^1$ branched along connected smooth algebraic curves
\cite{ADKY}, which includes an infinite family of surfaces of general type.

In all these examples it turns out that, if one considers projections of
sufficiently large degree (i.e., assuming $k\ge 3$ for $\CP^2$ and $k\ge 2$
for the other examples), the structure of $G^0$ is very simple, and obeys
the following conjecture:

\begin{conj}
Assume that $X$ is a simply connected algebraic surface and $k\gg 0$.
Then: $(1)$ the symplectic stabilization operation is trivial,
i.e.\ $K=\{1\}$ and $\bar{G}=\pi_1(\CP^2-D)$;
$(2)$ the homomorphism $\phi:\mathrm{Ab}(G^0)\to (\Z^2/\Lambda)^{N-1}$ is an
isomorphism; and
$(3)$ the commutator subgroup $[G^0,G^0]$ is a quotient of $\,\Z_2\times\Z_2$.
\end{conj}

\section{The symplectic isotopy problem} \label{sec:isotopy}

The symplectic isotopy problem asks under which conditions
(assumptions on degree, genus, types and numbers of singular points) it is
true that any symplectic curve in $\CP^2$ (or more generally in a complex
surface) is symplectically isotopic to a complex curve (by isotopy, we mean
a continuous family of symplectic curves with the same singularities).

The first result in this direction
is due to Gromov, who proved that every smooth symplectic curve of degree 1
or 2 in $\CP^2$ is isotopic to a complex curve \cite{Gr}. The argument
relies on a careful study of the deformation problem for pseudo-holomorphic
curves: starting from an almost-complex structure $J$ for which the given
curve $C$ is pseudo-holomorphic, and considering a family of
almost-complex structures $(J_t)_{t\in [0,1]}$ interpolating between $J$
and the standard complex structure, one can prove the existence of smooth
$J_t$-holomorphic curves $C_t$ realizing an isotopy between $C$
and a complex curve.

The isotopy property is expected to hold for smooth and nodal
curves in all degrees, and also for curves with sufficiently few cusps.
For smooth curves, successive improvements of Gromov's result have been
obtained by Sikorav (for degree $3$), Shevchishin
(for degree $\le 6$), and more recently Siebert and Tian \cite{ST}:

\begin{theorem}[Siebert-Tian]
Every smooth symplectic curve of degree $\le 17$ in $\CP^2$ is
symplectically isotopic to a complex curve.
\end{theorem}

Some results have been obtained by Barraud and Shevchishin for nodal
curves of low genus. For example, the following result holds \cite{Sh}:

\begin{theorem}[Shevchishin]
Every irreducible nodal symplectic curve of genus $g\le 4$ in $\CP^2$ 
is symplectically isotopic to a complex curve.
\end{theorem}

Moreover, work in progress by S.\ Francisco is expected to lead to an isotopy
result for curves of low degree with node and cusp singularities (subject
to specific constraints on the number of cusps).

If one aims to classify symplectic 4-manifolds by enumerating all branched
covers of $\CP^2$ according to the degree and number of singularities of
the branch curve, then the above cases are those for which the
classification is the simplest and does not include any non-K\"ahler
examples. On the other hand, Corollary \ref{cor:au2} implies that the
isotopy property cannot hold for all curves with node and cusp
singularities; in fact, explicit counterexamples have been constructed
by Moishezon \cite{Mo3} (see below).
\medskip

Even when the isotopy property fails, the classification of singular plane
curves becomes much simpler if one considers an equivalence relation weaker
than isotopy, such as {\it regular homotopy}, or {\it stable isotopy}.
Namely,
let $D_1,D_2$ be two Hurwitz curves (see Definition \ref{def:hurwitz}) in
$\CP^2$ (or more generally in a rational ruled surface), with node and
cusp singularities (or more generally singularities of type $A_n$).
Assume that $D_1$ and $D_2$ represent the same homology class, and that they
have the same numbers of singular points of each type. Then we have
the following results \cite{AKS,KK}:

\begin{theorem}[A.-Kulikov-Shevchishin]\label{thm:aks}
Under the above assumptions, $D_1$ and $D_2$ are {\em regular homotopic}
among Hurwitz curves, i.e.\ they are isotopic up to creations and
cancellations of pairs of nodes.
\end{theorem}

\begin{theorem}[Kharlamov-Kulikov]\label{thm:kk}
Under the above assumptions, let $D'_i$ $(i\in\{1,2\})$ be the curve
obtained by adding to $D_i$ a union of $n$ generic
lines (or fibers of the ruling) intersecting $D_i$ transversely at smooth
points, and smoothing out all the resulting intersections. Then for all
large enough values of $n$ the Hurwitz curves $D'_1$ and $D'_2$ are
isotopic.
\end{theorem}


Unfortunately, Theorem \ref{thm:aks} does not seem to have any implications
for the topology of symplectic 4-manifolds, because the node creation
operations appearing in the regular homotopy need not be admissible: even
if both $D_1$ and $D_2$ are branch curves of symplectic coverings, the
homotopy may involve plane curves for which the branched cover is not
smooth. For similar reasons, the applicability of Theorem \ref{thm:kk} to
branch curves is limited to the case of double covers, i.e.\ symplectic
4-manifolds which admit {\it hyperelliptic} Lefschetz fibrations. In
particular, for genus 2 Lefschetz fibrations we have the following result
\cite{AuGo}:

\begin{theorem}
If the symplectic 4-manifold $X$ admits a genus $2$ Lefschetz fibration,
then $X$ becomes complex projective after stabilization by fiber sums with
rational surfaces along genus $2$ curves.
\end{theorem}

It follows from Theorem \ref{thm:kk} that this result extends to all
Lefschetz fibrations with monodromy
contained in the hyperelliptic mapping class group.
However, few symplectic 4-manifolds admit such fibrations, and in general
the following question remains open:

\begin{question}
Let $X_1,X_2$ be two integral compact symplectic 4-manifolds with the same
$(c_1^2,\,c_2,\,c_1\!\cdot\![\omega],\,[\omega]^2)$. Do $X_1$ and $X_2$ become
symplectomorphic after sufficiently many fiber sums with the
same complex projective surfaces (chosen among a finite collection of model
holomorphic fibrations)?
\end{question}

This question can be thought of as the symplectic analogue of the classical
result of Wall which asserts that any two simply connected
smooth 4-manifolds with the same intersection form become diffeomorphic
after repeatedly performing connected sums with $S^2\times S^2$ \cite{Wall}.
\medskip

A closer look at the known examples of non-isotopic singular plane
curves suggests that an even stronger statement might hold.

It was first observed in 1999 by Fintushel and Stern \cite{FS2} that
many symplectic 4-manifolds contain
infinite families of non-isotopic smooth connected symplectic curves
representing the same homology class (see also \cite{Sm}). 
The simplest examples are obtained by ``braiding'' parallel copies of the
fiber in an elliptic surface, and
are distinguished by comparing the Seiberg-Witten invariants of
the corresponding double branched covers.
Other examples have been constructed by Smith, Etg\"u and
Park, and Vidussi. However, for singular plane curves the first examples
were obtained by Moishezon more than ten years ago \cite{Mo3}:

\begin{theorem}[Moishezon]
For all $p\ge 2$, there exist infinitely many pairwise non-isotopic
singular symplectic curves of degree $9p(p-1)$ in $\CP^2$ with
$27(p-1)(4p-5)$ cusps and $\frac{27}{2}(p-1)(p-2)(3p^2+3p-8)$
nodes, not isotopic to any complex curve.
\end{theorem}

Moishezon's approach is purely algebraic (using braid monodromy
factorizations), and very technical; the curves that he constructs
are distinguished by the fundamental groups of their complements
\cite{Mo3}. However a much simpler geometric description of this
construction can be given in terms of braiding operations, which makes
it possible to distinguish the
curves just by comparing the canonical classes of the associated branched
covers \cite{ADK}.

\begin{figure}[t]
\centering 
\epsfig{file=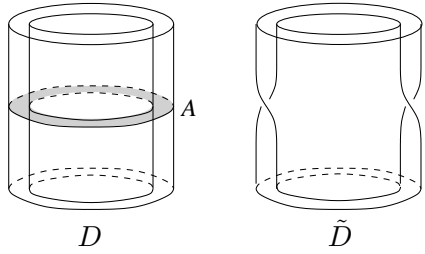,height=2.8cm}\\
$D$\hskip3cm$\tilde{D}$\vskip-2mm
\caption{The braiding construction}\label{fig:braiding}
\end{figure}

Given a symplectic covering $f:X\to Y$ with branch curve $D$, and given
a Lagrangian annulus $A$ with interior in $Y\setminus D$ and boundary contained in $D$,
we can {\it braid} the curve $D$ along the annulus $A$ by performing the
local operation depicted on Figure \ref{fig:braiding}. Namely, we cut out
a neighborhood $U$ of $A$, and glue it back via a non-trivial diffeomorphism
which interchanges two of the connected components of $D\cap \partial U$,
in such a way that the product of $S^1$ with the trivial braid is replaced
by the product of $S^1$ with a half-twist (see \cite{ADK} for details).

Braiding the curve $D$ along the Lagrangian annulus $A$ affects the branched cover
$X$ by a {\it Luttinger surgery} along a smooth embedded Lagrangian torus
$T$ which is one of the connected components of $f^{-1}(A)$ \cite{ADK}. This operation consists of cutting out 
from $X$ a tubular neighborhood of $T$, foliated by parallel
Lagrangian tori, and gluing it back via a symplectomorphism wrapping the
meridian around the torus (in the direction of the preimage of an arc joining
the two boundaries of $A$), while the longitudes are not affected.

The starting point of Moishezon's construction is the complex curve $D_0$
obtained by considering $3p(p-1)$ smooth cubics in a pencil, removing balls
around the 9 points where these cubics intersect, and inserting into each
location the branch curve of a generic degree $p$ polynomial map from
$\CP^2$ to itself. By repeatedly braiding $D_0$ along a well-chosen Lagrangian
annulus, one obtains symplectic curves $D_j$, $j\in\Z$.
Moishezon's calculations show that, whereas for the initial curve the
fundamental group of the complement $\pi_1(\CP^2-D_0)$ is infinite, the
groups $\pi_1(\CP^2-D_j)$ are finite for all $j\ne 0$, and of different
orders \cite{Mo3}. On the other hand, it is fairly easy to check that,
as expected from Theorem \ref{thm:adky}, this change in fundamental groups
can be detected by considering the canonical class of the
$p^2$-fold covering $X_j$ of $\CP^2$ branched along $D_j$. Namely, the
canonical class of $X_0$ is proportional to the cohomology class of the
symplectic form induced by the branched covering:
$c_1(K_{X_0})=\lambda[\omega_{X_0}]$, where $\lambda=\frac{6p-9}{p}$. On the
other hand, $c_1(K_{X_j})=\lambda[\omega_{X_j}]+\mu\,j\,[T]^{PD}$, where
$\mu=\frac{2p-3}{p}\neq 0$, and the homology class $[T]$ of the Lagrangian torus
$T$ is not a torsion element in $H_2(X_j,\Z)$~\cite{ADK}.
\medskip

Many constructions of non-K\"ahler symplectic 4-manifolds can be thought of
in terms of twisted fiber sum operations, or Fintushel-Stern surgery along
fibered links. However the key component in each of these constructions
can be understood as a particular instance of Luttinger surgery;
so it makes sense to ask to what extent Luttinger surgery may be responsible
for the greater variety of symplectic 4-manifolds compared to complex
surfaces. More precisely, we may ask the following questions:

\begin{question}
Let $D_1,D_2$ be two symplectic curves with nodes
and cusps in $\CP^2$, of the same degree and with the same numbers of
nodes and cusps.
Is it always possible to obtain $D_2$ from $D_1$ by a sequence
of braiding operations along Lagrangian annuli?
\end{question}

\begin{question}
Let $X_1,X_2$ be two integral compact symplectic 4-manifolds with the same
$(c_1^2,\,c_2,\,c_1\!\cdot\![\omega],\,[\omega]^2)$. Is it always possible
to obtain $X_2$ from $X_1$ by a sequence of Luttinger surgeries?
\end{question}

This question is the symplectic analogue of a question asked by Ron Stern
about smooth 4-manifolds, namely whether any two simply connected smooth 4-manifolds
with the same Euler characteristic and signature differ from each other by a
sequence of logarithmic transformations. However, here we do not require
the manifolds to be simply connected, we do not even require them to have
the same fundamental group.

\end{document}